\documentclass{gtart}

\input gtoutput
\volumenumber{3}\papernumber{6}\volumeyear{1999}
\pagenumbers{137}{153}\published{20 June 1999}
\proposed{David Gabai}\seconded{Walter Neumann, Cameron Gordon}
\received{1 September 1998}\revised{9 April 1999}
\accepted{14 June 1999}

\usepackage{amsmath,amssymb,amscd}

\newtheorem{thm}{Theorem}[section]
\newtheorem{lem}[thm]{Lemma}

\theoremstyle{definition}
\newtheorem{defn}{Definition}[section]

\theoremstyle{remark}
\newtheorem{rmk}{Remark}[section]

\def\R{\Bbb R}
\def\Z{\Bbb Z}

\def\F{\mathcal{F}}
\def\G{\mathcal{G}}

\newenvironment{pf}{\begin{proof}}{\end{proof}}
\let\Bbb\mathbb

\begin{document}

\title{$\R$--covered foliations of hyperbolic $3$--manifolds}
\asciititle{R-covered foliations of hyperbolic 3-manifolds}
\covertitle{R--covered foliations of hyperbolic 3--manifolds}

\author{Danny Calegari}
\address{Department of Mathematics \\ UC Berkeley \\ Berkeley, CA 94720}
\email{dannyc@math.berkeley.edu}
\primaryclass{57M50, 57R30}
\secondaryclass{53C12}
\keywords{$\R$--covered foliations, slitherings, 
hyperbolic $3$--manifolds, transverse geometry}
\asciikeywords{R-covered foliations, slitherings, 
hyperbolic 3-manifolds, transverse geometry}

\begin{abstract}
We produce examples of taut foliations of hyperbolic $3$--manifolds
which are $\R$--covered but not uniform --- 
ie the leaf space of the universal cover is $\R$, but pairs of leaves
are not contained in bounded neighborhoods of each other. This answers in the
negative a conjecture of Thurston in \cite{wT97}. We further show that
these foliations can be chosen to be $C^0$ close to foliations
by closed surfaces. Our construction underscores the importance of the
existence of transverse regulating vector fields and cone fields for
$\R$--covered foliations. Finally, we discuss the effect of perturbing
arbitrary $\R$--covered foliations.
\end{abstract}
\asciiabstract{%
We produce examples of taut foliations of hyperbolic 3-manifolds which
are R-covered but not uniform --- ie the leaf space of the universal
cover is R, but pairs of leaves are not contained in bounded
neighborhoods of each other. This answers in the negative a conjecture
of Thurston `Three-manifolds, foliations and circles I'
(math.GT/9712268). We further show that these foliations can be chosen
to be C^0 close to foliations by closed surfaces. Our construction
underscores the importance of the existence of transverse regulating
vector fields and cone fields for R-covered foliations. Finally, we
discuss the effect of perturbing arbitrary R-covered foliations.}

\maketitlepage

\section{$\R$--covered foliations that are not uniform}

\subsection{$\R$--covered foliations}

\begin{defn}
We say a foliation $\F$ of a compact $3$--manifold $M$ is {\em $\R$--covered}
if the pulled back foliation $\tilde \F$ of the universal cover
$\tilde M$ of $M$ is topologically the standard foliation of $\R^3$ by
horizontal $\R^2$'s.
\end{defn}

The first step in our construction is to produce a manifold $M$ with an
$\R$--covered foliation $\F$ which is not uniform. The condition that a
foliation be $\R$--covered is a somewhat elusive one, and in particular
it does not seem to be sufficient to find a cover $\hat M$ of $M$
so that the leaf space of $\hat \F$ is $\R$. This is related to the
question of when an infinite braid in $\R^3$ is trivial.

\begin{lem}
A taut foliation $\F$ of $M$ is $\R$--covered iff $\F$ has no spherical
or projective plane leaves, and for every arc $\alpha$ between
two points $p,q \in M$ there is an arc $\hat \alpha$ homotopic to
$\alpha$ rel.\ endpoints which is either contained in a leaf of $\F$ or
which is everywhere transverse to $\F$.
\end{lem}
\begin{pf}
If $\F$ is $\R$--covered, lift $\alpha$ to some $\tilde \alpha$ in $\tilde M$
and make it transverse there. If $\F$ is not $\R$--covered, either
$M$ is covered by $S^2 \times S^1$ or the leaf
space $L$ of $\tilde \F$ is a non-Hausdorff simply connected $1$--manifold.
This follows from a theorem of Palmeira, which says that the foliation
of $\R^3$ by the universal cover of a taut foliation is topologically
equivalent to a product of a foliation on $\R^2$ and $\R$, and therefore
such a foliation of $\R^3$ by planes whose leaf space is 
$\R$ is a product foliation (see \cite{cP78}).
In particular, there are distinct leaves $\lambda,\mu \in L$ which cannot 
be joined by an immersed path in $L$. That is, the topology of $L$ provides
an obstruction to finding such a $\hat \alpha$ as above.
\end{pf}

Note that one does not have an analogue of Palmeira's theorem for arbitrary
open $3$--manifolds --- that a foliation with leaf space $\R$ should be a
product --- and in fact this is not true. For example, remove from
$\R^3$ foliated by horizontal $\R^2$'s a properly embedded
bi-infinite transverse curve which does not intersect every leaf. This
is where the difficulty resides in showing that a foliation is $\R$--covered
by investigating an intermediate cover.

Many (most?) taut foliations of $3$--manifolds are not $\R$--covered. In
particular, by a theorem of Fenley, an $\R$--covered foliation of a
hyperbolic $3$--manifold has the property that in the universal cover,
every leaf limits to all of $S^2_\infty$. However, any compact leaf of a
taut foliation which is not a fiber of a fibration over $S^1$ lifts to a
quasi-isometrically embedded plane in the universal cover, 
and its limit set is a quasicircle (see \cite{sF91} for a fuller discussion).

\subsection{Uniform foliations}

\begin{defn}
A taut foliation $\F$ of a compact $3$--manifold $M$ is {\em uniform} if in the
pullback foliation $\tilde \F$ of the universal cover $\tilde M$, every
two leaves
$\lambda,\mu$ are a bounded distance apart. That is, there is
some $\epsilon$ depending on $\lambda,\mu$ so that $\lambda$ is contained in
the $\epsilon$--neighborhood of $\mu$, and vice versa. A foliation
$\F$ is obtained from a {\em slithering over $S^1$} if there is a 
fibration $\phi\co \tilde M \to S^1$ such that $\pi_1(M)$ acts as bundle maps
of this fibration, and such that the foliation of $\tilde M$ by components of
the fibers of $\phi$ agrees with $\tilde \F$. We will also refer to such
a foliation, perhaps ungrammatically, as a slithering.
\end{defn}

For additional details and definitions, see \cite{wT97}.
It is shown in \cite{wT97} that a uniform foliation with every leaf dense
is obtained from a slithering.

It is almost tautological from the definition of a slithering over $S^1$
that the action of $\pi_1(M)$ on the leaf space $L$ of $\tilde \F$ is
conjugate to a representation in $\widetilde {Homeo(S^1)}$, the
universal central extension of $Homeo(S^1)$. 

$3$--manifold topologists will be familiar with the short exact sequence
$$0 \to \Z \to \widetilde{PSL(2,\R)} \to PSL(2,\R) \to 0$$
This sits inside the short exact sequence
$$0 \to \Z \to \widetilde{Homeo(S^1)} \to Homeo(S^1) \to 0$$
Informally, $\widetilde{Homeo(S^1)}$ is the group of homeomorphisms of
$\R$ which are periodic with period $1$.

Let $Z$ be the generator of the center of $\widetilde {Homeo(S^1)}$. Then
$Z$ acts on $L$ by translations, and by choosing an invariant metric on
$L$ for this action, the action of every element of $\pi_1(M)$ on $L$
is periodic with some period equal to the translation length of $Z$. $Z$
is known as {\em the slithering map}.

The condition that a foliation be uniform is reflected in
the action of $\pi_1(M)$ on $L$ in the following way: since leaves
do not converge or diverge too much at infinity, holonomy cannot 
expand or contract the leaf space too much. 
Every compact interval in $L$
can be included in a larger compact interval which is ``incompressible'':
no element of $\pi_1(M)$ takes it to a proper subset or superset of itself.
If every leaf is dense, $L$ can be tiled with a countable collection of these
incompressible intervals laid end to end, so that the slithering map $Z$
acts on this tiling as a permutation.

Note that a foliation $\F$ may come from a slithering in many different
ways. For instance, if $\F$ admits a nonsingular transverse measure 
then the leaf space
of $\tilde \F$ inherits an invariant measure making it isometric to
$\R$. Then for any real 
$t > 0$ there is a slithering
$$\phi_t\co  \tilde M \to S^1$$ defined by the composition
$$\tilde M \to L \to S^1$$
where the last map is reduction mod $t$. Since the action of $\pi_1(M)$ on
$L$ preserves the property of points being integral multiples of $t$ 
apart, this map is a slithering.

However {\em generic} foliations come from a slithering in essentially
at most one way. The slithering
map $Z$ commutes with the action of every element of $\pi_1(M)$. So, for
instance, if an element $\alpha$ acts on $L$ with isolated fixed points,
the map $Z$ must permute this fixed point set. If, further, the action of
$\pi_1(M)$ on $L$ is minimal (ie every leaf of $\F$ is dense), the map
$Z$ is determined uniquely up to taking iterates. That is, there is
a {\em minimal} slithering $\phi\co \tilde M \to S^1$ with the property that
for every other slithering $\phi'\co \tilde M \to S^1$ determining the same
foliation, there is a finite cover $\psi\co S^1 \to S^1$ so that
$\phi = \psi \circ \phi'$.

This theory is developed in \cite{wT97}.

Following \cite{wT97} we define some auxiliary structure that will be used
to show that certain foliations are uniform or $\R$--covered.

\begin{defn}
Let $X$ transverse to $\F$ be a vector field. Then $X$ is {\em regulating}
if the lifts of the integral curves of $X$ to $\tilde M$ intersect every
leaf of $\tilde \F$.
\end{defn}

These lifts determine a one dimensional foliation of $\tilde M$. A leaf
in this foliation and a leaf in $\tilde \F$ intersect in exactly one point,
and consequently one can identify the leaf space of the one dimensional
foliation with any fixed leaf of $\tilde \F$---that is, with $\R^2$. Such
one dimensional foliations are called {\em product covered} in
\cite{CL98}. The main point for our purposes of this structure is the
following theorem:

\begin{thm}\label{excise_loop}
Suppose $\F$ is uniform (respectively $\R$--covered) and $X$ is a transverse
regulating vector field with a closed trajectory $\alpha$. Then the
restriction of $\F$ to $M - \alpha$ is also uniform (respectively $\R$--covered).
\end{thm}
\begin{pf}
The universal cover $\tilde M$ of $M$ is foliated as a product by
$\tilde \F$ and the integral curves of $X$ give this the structure of
a product $\R^2 \times \R$ in such a way that $\pi_1(M)$ acts by
elements of $Homeo(\R^2) \times Homeo(\R)$. Let $N$ be obtained from
$\tilde M$ by removing the lifts of $\alpha$. Then $N$ is the cover of
$M - \alpha$ corresponding to the subgroup of $\pi_1(M - \alpha)$
normally generated by the meridian of $\alpha$. One sees from the
structure of $X$ that $N$ is foliated as a product by infinitely punctured
disks and therefore that $\tilde N$ is $\R^3$ foliated by horizontal
$\R^2$'s. If $\F$ was uniform, any two leaves in $N$ would be a finite
distance apart. But leaves in $\tilde N$ correspond bijectively with
leaves in $N$ under the covering projection and therefore the same is
true in $\tilde N$; that is, the restriction of $\F$ to $M - \alpha$ is
uniform.
\end{pf}

Notice from the construction that if $\F$ came from a slithering, then
the restriction of $\F$ to $M - \alpha$ comes from a slithering 
which agrees with the restriction of the slithering map on $\tilde M$ to
the complement of the lifts of $\alpha$.
  
\subsection{Building uniform foliations from representations}

Let $F_g^n$ denote the surface of genus $g$ with $n$ punctures. Then
$F_1^1$ is the punctured torus, and $\pi_1(F_1^1) = \Z * \Z$. Let
$\alpha_l,\beta_l$ be standard generators for $\pi_1(F_1^1)$. Then we can
choose a representation $\rho$ of $\pi_1(F_1^1) \to \widetilde {Homeo(S^1)}$ by
sending $\alpha_l$ to translation through length $t$ and $\beta_l$ to 
some monotone element perhaps with a periodic collection of fixed points, each
distance $1$ apart. Let $M_l$ be the trivial circle bundle 
$M_l = F_1^1 \times S^1$ over the punctured torus, and pick
a flat $Homeo(S^1)$ connection on this bundle whose holonomy realizes the
representation $\rho$. Note that after fixing a trivialization of the
product, the representation is well-defined in $\widetilde {Homeo(S^1)}$
and not just $Homeo(S^1)$. The distribution determined by this connection
is integrable, by flatness, and integrates to give a foliation $\F_l$.

Geometrically, there is a foliation of $\tilde F_1^1 \times S^1$ by
leaves $\tilde F_1^1 \times {\text point}$. $\pi_1(F_1^1)$ acts on this
space by 
$$(x,\theta) \to (\alpha(x),\rho(\alpha)(\theta))$$ which preserves
the foliation. This foliation therefore descends to a foliation on
$$M_l = \tilde F_1^1 \times S^1 / \pi_1(M)$$
transverse to the $S^1$ fibers.

This foliation of $M_l$ comes from an obvious slithering
$\phi\co \tilde F_1^1 \times \R \to S^1$ which is just projection onto the
second factor followed by the covering map $\R \to S^1$ corresponding to
the circles in $M_l$.
The action of $\pi_1(M_l)$ on the leaf space is exactly
given by the representation $\rho$ in $\widetilde{ Homeo(S^1)}$ thought
of as sitting in $Homeo(\R)$.

Then the curve $\alpha_l \times 0$ sits in $M_l$ transverse to $\F_l$, 
and there is also a foliation on $M_l - \alpha_l$ which we also denote 
by $\F_l$.

It is easy to see that for irrational choice of $t$ the foliation $\F_l$ 
has every leaf dense. 

Furthermore we have the following lemma:

\begin{lem}\label{regulates}
The foliation $\F_l$ of $M_l - \alpha_l$ comes from a slithering. Furthermore,
this slithering can be taken to be the restriction of $\phi\co \tilde M_l \to S^1$
to the complement of the lifts of $\alpha_l$ in $\tilde M_l$. Moreover,
by choosing $\rho(\beta_l)$ suitably generic, this slithering is minimal
as defined above. 
\end{lem}
\begin{pf}
By theorem~\ref{excise_loop} it suffices to show 
there is a regulating vector field
of $M_l$ which agrees with $\alpha_l'$ when restricted to $\alpha$. Since
$M_l$ is topologically just $F_1^1 \times S^1$ it has a projection map to
$F_1^1$. Let $H$ be the torus foliated by circles that is the preimage of
the curve $\alpha$ on $F_1^1$ under this projection. Then $M_l - N(H)$ has
an obvious codimension $2$ foliation by $S^1$ fibers. We extend this
foliation over $N(H)$, which can be parameterized as 
$S^1 \times S^1 \times [-1,1]$, by foliating each $S^1 \times S^1 \times *$
with parallel lines which rotate continuously from vertical (parallel to
the $* \times S^1$ direction in $M_l$) on the boundary to horizontal
(parallel to $\alpha$) on $H$, always staying transverse to $\F_l$. It is
obvious that this is a foliation by regulating curves, and we denote
its associated unit tangent vector field by $X_l$.

If we choose $\rho(\beta_l)$ to be generic and
close to the identity with isolated fixed points, the slithering defined
in the statement of the theorem is minimal.
\end{pf}

On another punctured torus with basis for $\pi_1$ given by
$\alpha_r,\beta_r$ we pick
another representation $\sigma$ in $\widetilde{Homeo(S^1)}$ 
so that $\sigma(\alpha_r)$ is translation
through $s$, where again $s$ is irrational and incommensurable with $t$,
and $\sigma(\beta_r)$ is some random element which commutes with $Z$ but not
with $\alpha_r$. Then we can form $M_r = F_1^1 \times S^1$ foliated as
above, and remove $\alpha_r$ from $M_r$ to produce another foliated manifold
with a slithering.

Let $M$ be obtained by gluing $M_l-\alpha_l$ and $M_r-\alpha_r$ along the
torus boundaries of neighborhoods of $\alpha_l$ and $\alpha_r$ respectively.
Denote this torus in the sequel by $S \subset M$. Each piece $M_l - \alpha_l$,
$M_r-\alpha_r$ admits a regulating vector field $X_l,X_r$ as constructed in
lemma~\ref{regulates}.
We perform this gluing in such a way that the foliations of the boundary
tori by meridional circles agree. If we like, we can perform the
gluing so that the restriction of the slithering maps for the left and
the right foliations commute, when restricted to their action on the
leaf space of the universal cover of $S$.

Then $$\pi_1(M) = \pi_1(M_l-\alpha_l) *_{\Z \oplus \Z} \pi_1(M_r-\alpha_r)$$
acts on $\R$ by the amalgamated action of each piece on the leaf space
of its respective universal cover, each canonically identified with the
leaf space of the foliation of $\tilde S$.

Topologically, $M$ is a graph manifold obtained from four copies of
$F_0^3 \times S^1$. To see this, observe that $M_l - N(H)$ (with notation
as above) is exactly $(F_1^1 - N(\alpha)) \times S^1$ 
which is $F_0^3 \times S^1$.
Also, observe that $N(H)= S^1 \times S^1 \times [-1,1]$ in many different ways,
including a way in which $\alpha$ is $S^1 \times * \times 0$.
It follows that $N(H) - \alpha$ is topologically also $F_0^3 \times S^1$.
However, these foliations by circles cannot be made to agree on the
boundary tori of different pieces, and $M$ is not a Seifert fibered space.

Let $\F$ denote the induced foliation of $M$. Is $\F$ $\R$--covered?
To establish that it is indeed $\R$--covered, it will suffice to show that each
piece is uniform and admits a regulating transverse vector field $X_l,X_r$
which agree on the gluing torus to make a global regulating vector field $X$.

\begin{lem}
$\F$ is an $\R$--covered foliation of $M$.
\end{lem}
\begin{pf}
Let $\hat M$ foliated by $\hat \F$ 
be the cover of $M$ obtained by taking copies of the universal
covers of $M_l$ and $M_r$, drilling out countably many copies of the lifts of
$\alpha_l$ and $\alpha_r$, then gluing along the boundary components.
Then the regulating vector fields $X_l,X_r$ lift to $\hat M$ to give a
global trivialization of this manifold as a product of an infinite genus
Riemann surface with $\R$. This implies that the universal cover of
$\hat M$ is the standard $\R^3$ foliated by $\R^2$'s, and we see therefore
that $\F$ is $\R$--covered.
\end{pf}

Since every leaf of $\F$ is dense in $M$, if $\F$ were uniform it would come 
from a slithering by \cite{wT97}. 
However, we have seen that the action of $\pi_1(M)$ on the
leaf space of the universal cover is the amalgamation of the actions
of $\pi_1(M_l - \alpha_l)$ and $\pi_1(M_r - \alpha_r)$ along their
gluing $\Z \oplus \Z$. It follows that there is no single translation 
$\mu \in Homeo(\R)$ (ie an element without fixed points) 
which commutes with both $\beta_l$ and $\beta_r$ for sufficiently generic
choice of $\sigma(\beta_r)$ and $\rho(\beta_l)$, since the periods of
the left and right slithering maps are incommensurable.

More explicitly, let us fix a lift $\tilde S$ 
of the torus $S$ to $\tilde M$ which divides
a piece of $\tilde M$ which is a lift of $M_l - \alpha_l$ from
a piece which is a lift of $M_r - \alpha_r$. Let us identify
the leaf space $L$ with the leaf space of $\tilde S$. Let $\alpha \in
\pi_1(M)$ corresponding to the longitude of $S$ preserve $\tilde S$ and
act on $L$ as a translation. Fix a metric on $L$ such that
$\alpha$ acts as translation through a unit length. 
Let $Z_l$ and $Z_r$ be the translations in $Homeo(L)$ 
corresponding to the slithering map of $M_l - \alpha_l$ and $M_r - \alpha_r$
thought of as acting on the leaf space of $\tilde S$. Then $Z_l$ acts
as translation through length $\frac 1 s$ and $Z_r$ acts as translation 
through length $\frac 1 t$. By minimality, the only translations in
$Homeo(L)$ that commute with $\rho(\beta_l)$ are multiples of $Z_l$,
whereas the only translations that commute with $\sigma(\beta_r)$ are
multiples of $Z_r$. It follows that no translation commutes with both
elements, and $\F$ does not come from a slithering.

We have therefore proved the following theorem:

\begin{thm}
$\F$ as above is $\R$--covered but not uniform.
\end{thm}

\section{Lorentz cone fields}

The following definition is from \cite{wT97}:

\begin{defn}
On any manifold $M$, a {\em (Lorentz) cone field} $C$ transverse to a 
codimension one foliation
$\F$ is the field of timelike vectors (ie with positive norm) 
for a (continuously varying) form on
$TM$ of signature $(n-1,1)$ such that $T\F_p$ are spacelike. A cone field is
{\em regulating} if every complete curve $X$ with $X' \in C$ is
regulating for an $\R$--covered $\F$.
\end{defn}

Regulating cone fields are discussed in \cite{wT97}, and shown to exist
for all uniform foliations. We show that the example constructed in the
previous section admits a regulating cone field.

Each piece $M_l - \alpha_l$, $M_r -\alpha_r$ admits a regulating cone field 
$C_l,C_r$ which is degenerate along the boundary torus, coming from the
restriction of the regulating cone fields on $M_l,M_r$ which are tangent
to $\alpha_l$ and $\alpha_r$. 
Let $C$ denote the cone field on $M$ which agrees with $C_l$ and $C_r$
away from a collar of the separating torus, and which near this separating
torus is sufficiently thin so that every curve which crosses this collar
must wind a distance at least $T$, as measured in periods of the
longitude, transverse to the foliation.

\begin{thm}
$C$ as defined above is a regulating cone field for $\F$.
\end{thm} 
\begin{pf}
We will do our calculations in $\hat M$, using the fact that $\hat M$
admits a regulating vector field coming from $X_l,X_r$ which preserves
each uniform piece. Observe that $\hat M$ has a decomposition into a
countable collection of covering spaces of $M_l-\alpha_l$, $M_r-\alpha_r$
which we denote $L_i,R_i$ for some particular choice of indices $i$.
Also label the separating cylinders, all of them lifts of the gluing
torus in $M$, as $A_i$ for some index $i$. Let $d_i$ denote the
metric on $L$, the leaf space of $\hat \F$, given by the transverse
measure on $\F|_S$ scaled so that the curves $\alpha_l = \alpha_r$ have
period $1$ and so that the left and right slithering maps act by
translations of this measure, as measured in the cylinder $A_i$.
Then we claim the following lemma:

\begin{lem}
If $A_i$ and $A_j$ are separated by $n$ cylinders $A_k$ then
for any two leaves $\lambda,\mu \in L$, 
$$|d_i(\lambda,\mu) - d_j(\lambda,\mu)| \le (n+1) \max(\frac 1 s,\frac 1 t)$$
\end{lem}
\begin{pf}
The proof follows immediately by induction once we show the result for
$n=0$. If $A_i$ and $A_j$ bound a single $L_i$ or $R_i$, then 
the fact that the pieces $L_i$ and $R_i$ slither over $S^1$ implies that 
$d_i(\lambda,\mu)$ and $d_j(\lambda,\mu)$ differ by at most one period of
the slithering, which in terms of the measure on the boundary torus, is
either $\frac 1 s$ or $\frac 1 t$ depending on whether we are in an $L_i$
or an $R_i$. 
\end{pf}

Let $r = \max (\frac 1 s,\frac 1 t)$.

Now we will show that any curve $\gamma$ 
supported by $C$ makes definite progress
relative to any given $d_i$, and therefore relative to all of them.
In particular, it is regulating. After re-ordering indices,
any such infinite curve in $\hat M$, starting on some leaf $\lambda \in L_1$,
can be broken up into segments $\gamma_1,\gamma_2,\dots$ where each
$\gamma_i$ is contained in $L_i$ or $R_i$ (according to sign). It is clear
that if there are only finitely many $\gamma_i$ (ie $\gamma$ crosses
only finitely many separating annuli) then eventually $\gamma$ can
be seen to be regulating, since it is supported by some lift of $C_r$ or
$C_l$. So we suppose there are infinitely many $\gamma_i$. Let $Y_i$ be
the union of the first $i$ segments of $\gamma$.
Let $Z_i$ be the {\em shadow} of
$Y_i$ on $\lambda$; that is, the curve on $\lambda$ obtained by projecting
$\hat M$ to $\lambda$ along the integral curves of $X$. Let $Y_i'$ be
the integral curve of $X$ interpolating between the endpoint of $Z_i$ and
the endpoint of $Y_i$. Then $d_i(Y_j) = d_i(Y_j')$ for all $i,j$ since
the curves $Y_j$ and $Y_j'$ have endpoints on the same pair of leaves.

We can estimate $d_1(Y_1') \ge T$ by hypothesis
on $C$. It follows that $d_2(Y_1') \ge T - r$. But then $d_2(Y_2') \ge
2T - r$ and so $d_3(Y_2') \ge 2T - 2r$. Continuing in this way and by
induction, we get $d_n(Y_n') \ge nT - nr$. But then 
$$d_1(Y_n) = d_1(Y_n') \ge nT - 2nr$$ 
and one can see that by choosing $T \gg 2r$ the curve $\gamma$ makes 
arbitrary progress relative to some {\em fixed} $d_i$, and therefore
is regulating.
\end{pf}

An instructive analogy is given by the comparison between imperial and
metric weights and measures: suppose I have a small object which I can
measure to the nearest inch or to the nearest centimeter. Then I get
estimates which vary greatly compared to the length of the object. If the
object is much bigger, the estimates are comparatively better. My
regulating curve above makes definite progress, even though its progress
is translated into inches, then centimeters, then inches, then centimeters
\dots rounding down every time. 

\subsection{A hyperbolic example}

\begin{thm} \label{T:lorentz}
Suppose $M$,$\F$ is any compact oriented $3$--manifold with a co-orientable
$\R$--covered foliation, and suppose that $\F$ admits a transverse regulating 
Lorentz cone field $C$. Let
$\gamma$ be any simple closed curve supported by $C$. If $M_n(\gamma)$
is obtained by taking an $n$--fold branched cover over $\gamma$, and
$\F_n(\gamma)$ denotes the pullback foliation, then $\F_n(\gamma)$ 
is $\R$--covered. Moreover, $\F_n(\gamma)$ is uniform iff $\F$ is. 
\end{thm}
\begin{pf}
The point of having a regulating cone field $C$ is that for any $\gamma$
supported as above, there is a regulating vector field $X$ of $M$ so that
$X|_\gamma = \gamma'$. This follows immediately from obstruction theory,
once one notices that sections of $C$ are contractible; eg use a
partition of unity.

Now in $\tilde M$, the universal cover of $M$, $\gamma$ lifts to a
collection of bi-infinite regulating curves, and $\tilde M_n(\gamma)$ is
the universal orbifold cover of $\tilde M$ where we declare that there are
order $n$ cone angles along the lifts of $\gamma$. Let $\tilde \F_n(\gamma)$
be the pullback foliation in that universal orbifold cover and let
$\lambda,\mu$ be two leaves there. They can be joined by some arc $\alpha$
in the complement of the lifts of the cone locus, which projects to an arc
$\pi(\alpha)$ in $\tilde M$. By homotoping $\pi(\alpha)$ rel.\ endpoints
along integral curves of $X$, we can make it transverse to $\tilde \F$ 
{\em without crossing any lift of $\gamma$}. Then this perturbed 
$\pi(\alpha)$ lifts to a perturbed $\alpha$ in $\tilde M_n(\gamma)$ 
transverse to $\tilde \F_n(\gamma)$, thereby demonstrating that
$\F_n(\gamma)$ is $\R$--covered, as required.

If $\F$ was not uniform, there would be a pair of leaves in $\tilde \F$ which
diverge at infinity. They lift to leaves with the same property in
$\tilde \F_n(\gamma)$. Alternatively, the uniformity or lack thereof can be
seen in the action of $\pi_1(M_n(\gamma))$ on $\R$.
\end{pf}

We return to the example $M$ that we constructed in the previous section.
In $M$, it is clear that we can choose a curve $\gamma$ supported by
$C$ whose complement is hyperbolic. For, we can certainly do this in each
side of $M$, and then by crossing back and forth across the separating
torus, we can arrange for the complement of $\gamma$ to be atoroidal. 
In particular, by choosing $\rho(\beta_l)$
and $\sigma(\beta_r)$ to be sufficiently close to translations, the
regulating cone fields $C_l$ and $C_r$ can be as ``squat'' as we like, and
we have a great deal of freedom in our choice of the restriction of 
$\gamma$ to the complement of a collar neighborhood of $S$.

For example, if $\rho, \sigma$ are chosen so that every element 
acts on the leaf space as a translation through a
rational distance, $\F$ would be a surface bundle over a circle; in a
product bundle, a curve whose projection to $S^1$ is a
homeomorphism and whose projection to the base surface fills up the
surface (ie complementary regions for the geodesic representative
are disks) has atoroidal complement. Similar curves are easily found in
any surface bundle, and one can arrange for them to wind several times
in the circle direction when they pass through some reducing torus.
The point is that any sufficiently complicated
curve will suffice. Then for nearby choices of $\rho,\sigma$, such a
curve will still be regulating and contained in the 
regulating Lorentz cone field, as we show in the following lemma:

\begin{lem}
Suppose $\F_l'$ is a transversely measured foliation of $M_l - \alpha_l$ as
above coming from some representation $\rho(\beta_l) = \text{translation}$, 
and let $C'$ be a transverse Lorentz cone field for $\F_l'$ appropriately
degenerate near $\alpha_l$. Then for
slitherings $\F_l$ coming from sufficiently close choices of $\rho(\beta_l)$,
the cone field $C'$ is regulating for $\F_l$.
\end{lem}
\begin{pf}
We need to check for sufficiently mild perturbations $\F_l$ of $\F_l'$ 
that any curve supported by $C'$ passes through at least one period of 
the slithering of $\F_l$, since then it must pass through 
arbitrarily many such periods and therefore be regulating.

Observe first that any transverse Lorentz cone field is regulating for
a transversely measured foliation, since one can uniformly compare
distance along a curve supported by the cone field and distance with
respect to the transverse measure.

For $\F_l$ sufficiently close to $\F_l'$, $C'$ is a transverse Lorentz
cone field for $\F_l$. The codimension $2$ foliation $X_l$ described 
in lemma~\ref{regulates} is regulating for every choice of 
$\rho(\beta_l)$, and we assume that this foliation lifts to the 
vertical foliation of $\R^3$ by $\text{point} \times \R$. We choose
co-ordinates on $\R^3$ so that the regulating curves are of the form
$x=\text{const.}$, $y=\text{const.}$ and the leaves of $\tilde \F_l'$ are
of the form $z = \text{const.}$

The slithering on $M_l - \alpha_l$ comes from the slithering on $M_l$
associated with the circle bundle as above. We assume that the
circles of the Seifert fibration away from $N(H)$ lift to vertical
arcs of length $1$.

With this
choice of co-ordinates, for any foliations $\F_l$ constructed from
a representation as above and for every $x,y,z$, the points $(x,y,z)$ and 
$(x,y,z+1)$ in the universal cover are on leaves which are one period of 
the slithering apart. 

For a point $p$ in $\tilde M_l$ given in co-ordinates by $(x,y,z)$, let
$\lambda_p$ be the leaf of $\tilde \F_l$ through $p$, and $\mu_p$ the leaf of
$\tilde \F_l$ through $(x,y,z+1)$. Then the leaves $\lambda_p$ and
$\mu_p$ differ by a translation parallel to the $z$-axis of length $1$.
Note that this translation need in no way correspond to the action of
an element of $\pi_1(M_l)$ on $\tilde M_l$. The light cone of
$C'$ through $p$ intersects the horizontal plane passing through $(x,y,z+1)$
in a compact region. For sufficiently small perturbations of
$\rho(\beta_l)$, the leaf $\mu_p$ will be a small perturbation of
the horizontal plane through $p$, and the intersection of the light
cone of $C'$ through $p$ with $\mu_p$ will also be a compact
region. 

Now, $M_l - \alpha_l$ is not compact, but we can consider its double
$N$, foliated by the double of $\F_l'$, and equipped with a transverse
Lorentz cone field obtained by doubling $C'$ which is degenerate along the
boundary components of $M_l - \alpha_l$. The compactness of $N$
implies that for a sufficiently small perturbation $\F_l'$ of $\F_l$, the
intersections as above will be compact for all $p$. 
This implies that a curve supported by $C'$ will need to go only a
bounded distance before traveling a full period of the slithering. In
particular, any bi-infinite curve supported by $C'$ will travel through
infinitely many periods of the slithering in either direction and will 
therefore be regulating, which is what we wanted to show.
\end{pf}

\begin{rmk}
What is really essential to notice in the above set up is that the leaves
passing through $(x,y,z)$ and $(x,y,z+1)$ in the universal cover were one 
period of the slithering away from each other for {\em all} $\F_l$. In some
sense, all the slitherings are determined by the structure of the Seifert
fibration where they originated. For a generic perturbation
of a uniform foliation, one has no control over how the slithering map
may vary, or even whether the perturbed foliations are uniform at all.
\end{rmk}

Since for a transversely measured foliation, {\em any} transverse Lorentz
cone field is regulating, we can choose our curve
$\gamma$ to be any curve transverse to a measured foliation $\F_l'$ with
hyperbolic complement which is sufficiently steep near the separating torus, 
and then choose the representation $\rho(\beta_l)$
to be sufficiently close to a translation that $\gamma$ stays in a
regulating cone field.

Since $M - \gamma$ is
hyperbolic, for sufficiently large $n$ an $n$--fold branched cover of
$M$ over $\gamma$ is hyperbolic. One cannot be sure that an $n$--fold
branched cover will always exist, but at least one has a hyperbolic
orbifold structure on $M$ with cone angle $2\pi/n$ along $\gamma$,
and by Selberg's lemma (see \cite{aS60}) 
one knows this is finitely (orbifold) covered by a
genuine hyperbolic manifold which is a branched cover of $M$ along
$\gamma$.

By the discussion above, the
induced foliation is $\R$--covered but not uniform. Moreover, by choosing
$\rho(\beta_l)$ and $\sigma(\beta_r)$ sufficiently close to translations,
the foliation $\F$ is as close to
a transversely measured foliation as we like. Passing to a branched cover
preserves this fact. Since transversely measured foliations of $3$--manifolds
are arbitrarily close (as $2$--plane fields)
to surface bundles over $S^1$, we have proved:

\begin{thm}
There exist foliations of hyperbolic $3$--manifolds which are $\R$--covered
but not uniform. They can be chosen arbitrarily close to surface bundles
over circles.
\end{thm}

\begin{rmk}
It is clear that the construction outlined above can be made in some
generality. One can construct $\R$--covered but not uniform foliations
by plumbing together uniform foliations along boundary tori in numerous ways.
In a great number of cases, these will admit regulating transverse cone fields,
and by branching as above one can produce many atoroidal examples.
One can easily arrange for these examples to be closed; for instance,
by doubling $M$ before removing a curve with atoroidal complement in the
examples constructed above.

Another construction, explained in detail in \cite{wT97}, involves choosing a
train track with integer weights supported by a regulating cone field, then
plumbing the leaves along the train track with a surface of genus given
by the track weight. These plumbed surfaces can be ``twisted'' by a
surface automorphism under the monodromy around loops of the train track. 
Thurston expects that these examples are
sufficiently flexible to allow one to prescribe the homology of $M$.
\end{rmk}

\section{General $\R$--covered foliations}

\subsection{Regulating vector fields}

For a general $\R$--covered foliation, we do not know whether or not
there exists a transverse Lorentz cone field.

In what follows, $\F$ will be an $\R$--covered taut foliation of a 
$3$--manifold $M$ with hyperbolic leaves.

We show in \cite{dC98} that the circles at infinity of each leaf in the
universal cover of $M$ can
be included in a topological {\em cylinder at infinity} $C_\infty$ on which
$\pi_1(M)$ acts by homeomorphisms. 

The following theorem is proved in \cite{dC98}: 
\begin{thm}
With notation as above, there is a
global trivialization of $C_\infty$ as $S^1 \times \R$ so that the
action of $\pi_1(M)$ preserves the horizontal and vertical foliations of
$C_\infty$ by $S^1 \times \text{point}$ and $\text {point} \times \R^1$.
\end{thm}

In \cite{wT98} it is suggested that {\em all} taut foliations should
have a pair of essential laminations $\Lambda_+, \Lambda_-$
transverse to each other and to
the foliation which intersect each leaf in a ($1$--dimensional) geodesic
lamination. In the case of $\R$--covered foliations, these laminations
should come from a pair of $1$--dimensional laminations 
$\hat \Lambda_+,\hat \Lambda_-$ of the universal
$S^1$ described in the theorem above which are invariant under the
action of $\pi_1(M)$. 

Let the universal $S^1$ bound a hyperbolic plane $D$, and
let $\hat \Lambda_+,\hat \Lambda_-$ be the associated geodesic laminations
in $D$. Then since each circle at infinity is canonically associated with
this $S^1$, each point in $\hat \Lambda_+ \cap \hat \Lambda_-$ determines
a unique point in each leaf of $\tilde \F$. Similarly, each segment of
$\hat \Lambda_\pm$ between points of intersection determines a unique
geodesic segment in each leaf, and each complementary region determines
a unique geodesic polygon in each leaf. If one fixes some canonical
geometric parameterization of the family of convex geodesic polygons of
a fixed combinatorial type, this parameterization gives rise to
a canonical identification of each leaf with each other leaf, preserving
the stratification outlined above. The fibers of this identification give
a one-dimensional foliation transverse to $\tilde \F$, and the tangent
vectors to this foliation are a regulating vector field.

It is easy to see that non-quadrilateral complementary regions in $D$
give rise to solid cylinders in $\tilde M$ which cover solid tori in
$M$, since their cores are isolated. 
The cores of these solid tori are necessarily regulating, and
therefore define elements of $\pi_1(M)$ which act on the leaf space by
translation.

In fact, Thurston has communicated to this author a sketch of a proof
that an $\R$--covered foliation admits {\em some} transverse 
lamination which intersects every leaf in geodesics (though not necessarily
a pair of such). In \cite{dC98} we can show by methods slightly
different to those above that this assumption is
enough to imply that there exists a regulating vector field transverse
to $\F$ which can be chosen to have closed orbits. In short, we have the
following theorem:

\begin{thm}
Let $\F$ be an $\R$--covered foliation. Then there exists a regulating
vector field transverse to $\F$ which can be taken to have closed orbits.
These orbits determine elements $\alpha \in \pi_1(M)$
which act on the leaf space of $\tilde \F$ without fixed points. Branched
covers of $M$ over these closed orbits give new $\R$--covered foliations.
\end{thm}

One hopes these results are all pieces of a unified picture tying the
intrinsic geometry of $\R$--covered foliations to the extrinsic geometry of
the foliated manifolds containing them. Not all the details of this picture
are yet visible. Nevertheless, it seems worthwhile to make this picture as 
explicit as possible.

\subsection{Instability of $\R$--covered foliations}

Despite the positive results of the previous sections, it seems that
the property of being $\R$--covered is quite delicate.
The following example is suggestive.

Let $M$ be a hyperbolic surface bundle over a circle with fiber $F$ and
pseudo-Anosov monodromy $\psi\co F \to F$. Let $\F$ be the induced foliation
by surfaces. Let $\gamma$ be a simple closed curve on $F$ so that
$\gamma \cap \psi(\gamma) = \emptyset$. Note that it is easy to show that
there exist such examples, by first choosing any $\gamma,\psi,F$ and then
using the fact that surface groups are LERF to lift to a finite regular
cover where $\gamma$ and its image are disjoint.

Let $\hat M$ be the $\Z$--cover of $M$ defined by the circle direction.
Topologically, $\hat M$ is $F \times \R$ foliated by closed surfaces 
$F \times \text{point}$.
Let the group $\Z$ generated by the deck translation, which we denote
$\Psi$, act by $\Psi(x,t) = (\psi(x),t+1)$.

Let $A$ be the annulus $\gamma \times [-\epsilon,1+ \epsilon]$. Then
by construction, $A$ and its translates are disjoint. Let 
$\mu\co [-\epsilon,1+ \epsilon] \to [-\epsilon,1 + \epsilon]$ be a
homeomorphism close to the identity which moves every point except the
endpoints up some small amount. Then we can cut open $\hat M$ along $A$
and its translates, and shear the foliation on one side by the 
translation $\text{id} \times \mu$ to get a new foliation $\G$. This
can certainly be done in such a way that $\G$ is arbitrarily close to $\F$.

Now whenever an integral curve of $\G$ passes through $A$ or a translate in
the ``positive'' direction, it will be sheared upward (relative to $\F$)
by $\mu$. 

Let us suppose that by choosing a suitable path in $\G$, we can arrange
that a curve starting at some $(x_0,0)$ can get to $(x_1,1)$ in length $t$, as
measured in $\G$, by
winding sufficiently many times through $A$. Now, the curve can continue
to wind around $\Psi(A)$, and after moving a distance $2t$, it can reach
$(x_2,2)$, and so on. Remember that there is a transverse regulating 
vector field $X$ to $\F$ given by tangents to the curves 
$\text{point} \times \R$.

Since $\psi$ is pseudo-Anosov, when we compare arclength at $(x,t)$ and
$(x,0)$ by projection along $X$, we see that a vector of norm $\lambda^t$
at $(x,t)$ might project to a vector of norm $1$ at $(x,0)$, where $\lambda$ is
the multiplier of $\psi$ on the invariant transverse measure of the
unstable lamination of $F$. Hence, as measured in $\G$, a curve
$\gamma$ beginning at $(x,0)$
could have length $nt$ but its projection to $F,0$ could have
length as little as $\sum_{i=1}^n t/{\lambda^i}$. In particular, a curve
in $\G$ could ``escape to infinity'' while its projection to $(F,0)$ could
move only a finite distance. By picking two points $(x,0)$ and $(x,n)$ 
sufficiently far apart, and moving them by curves in $\G$ joined 
    by integral curves of $X$, it seems plausible that one could find a 
path in $\G$ where holonomy was not defined after some finite time, 
suggesting that $\G$ was not $\R$--covered. 

Of course, there are problems with making this concrete: distances in
$\G$ are only magnified in the direction of the stable lamination as we
go upwards; perhaps to make the curves cross through $A$ and its translates
sufficiently often, we need to go in both stable and unstable directions.
Moreover, even if one could show that holonomy failed to be defined for
all time along integral curves of $X$, it does not rule out the possibility
that $\G$ is still $\R$--covered and $X$ is merely not regulating, 
although the author understands that very recently S\' ergio Fenley has shown
that a pseudo-Anosov flow transverse to an $\R$--covered foliation should
always be regulating (\cite{sF98}). A similar result is also proved in
\cite{dC98}.

\subsection*{Acknowledgments}

In writing this paper I benefited from numerous helpful conversations 
with Andrew Casson, S\' ergio Fenley and Bill Thurston. In particular,
many of the ideas contained here are either implicit or explicit in
the wonderful paper \cite{wT97}. While writing this paper, I was partially
supported by an NSF Graduate Fellowship.
\end{document}